\documentclass[12pt]{amsart}
\usepackage{amscd,amssymb}
\usepackage{latexsym}
\DeclareFontEncoding{OT2}{}{} 
  \newcommand{\textcyr}[1]{%
    {\fontencoding{OT2}\fontfamily{wncyr}\fontseries{m}\fontshape{n}%
     \selectfont #1}}
\newcommand{\Sha}{{\mbox{\textcyr{Sh}}}}

\newcommand{\bpf}{\noindent {\bf Proof. }}
\newcommand{\epf}{$\Box$ \vspace{+10pt}}
\newcommand{\rk}{\noindent {\bf Remark. }}

\newtheorem{theorem}{Theorem}[section]
\newtheorem{proposition}[theorem]{Proposition}
\newtheorem{lemma}[theorem]{Lemma}
\newtheorem{corollary}[theorem]{Corollary}

\newtheorem{definition}{Definition}[section]

\newcommand{\A}{\mathbb A}
\newcommand{\Q}{{\mathbb Q}}

\newcommand{\PP}{{\mathbb P}}

\newcommand{\N}{{\mathbb N}}
\newcommand{\F}{{\mathbb F}}

\newcommand{\Z}{{\mathbb Z}}
\newcommand{\V}{V_{\alpha,\beta}}
\newcommand{\C}{C_{\alpha,\beta}}

\newcommand{\LL}{{\mathcal{L}}}

\newcommand{\TT}{{\mathcal T}}

\newcommand{\GG}{{\mathcal G}}

\newcommand{\oh}{\mathcal O}
\newcommand{\Kbar}{\overline{K}}

\newcommand{\pgl}{\mbox{PGL}}
\newcommand{\gl}{\mbox{GL}}

\begin{document}

\pagestyle{plain}

\title{Explicit Descent over $X(3)$ and $X(5)$}

\author{Catherine O'Neil}

\date{\today}

\begin{abstract}

Let $E$ be an elliptic curve over a field $K,$ and let $p$ denote either 3 or
5. In the case that all of the $p$-torsion points of $E$ are defined over $K,$
we give an explicit description of $p$-descent on $E$ in two parts.  The first
part consists of an explicit map from the finite abelian group $E(K)/p E(K)$ to
$Sel(E)[p],$ and the second part is an explicit identification of $Sel(E)[p]$
with models of principal homogeneous spaces of $E$ inside $\PP^{\, p-1}.$

\end{abstract}

\maketitle

\section{Introduction}

This paper represents a complete solution, in a special case, to the problem of
performing explicit $p$-descent on elliptic curves for the primes $p=3$ and
$p=5.$  Namely, we work with the assumption that all of the $p$-torsion points
of the elliptic curve $E$ are defined over the field of definition of $E.$ 
From this vantage point, an obstruction (called the `period-index obstruction')
to constructing smooth models of genus one curves in low-dimensional projective
spaces is extremely concrete and in some sense is the only obstruction to
performing descent.   The general case, when there are no rationality
assumptions on the $p$-torsion, also contains this obstruction and much more;
we believe it is instructive to isolate this obstruction, i.e. to consider this
a kind of `base case' for a general method.

A few remarks are warranted.  First, 2-descent over $\Q$ is completely known
and implemented for the computer with no rationality assumptions on the
2-torsion points.  References for 2-descent include John Cremona's `mwrank'
program and \cite{cremona} (see pages 84 and 85 for the manifestation of the
period-index obstruction in the case of 2-descent), based on the methods
explained in \cite{Birch}.  For complete 2-descent over a general number  field
see \cite{simon1}.  An algorithmic method for $p$-descent with no rationality
assumption and which generalises to some higher genus cases can be found in
\cite{stolletal}.

Next, why do 3-descent if 2-descent is known?  Loosely speaking, when there is
nontrivial 2-torsion in $\Sha,$ 2-descent is not effective, but 3-descent may
be, if the 3-torsion of $\Sha$ is trivial.   In other words the obstruction to
the effectiveness of $p$-descent is $\Sha[p].$  Since $\Sha$ is conjecturally
finite, descent will be `eventually effective.' Certainly it is useful to have
two or three primes completely understood, and maybe more if possible.

Next, we know (by the Galois equivariance of the Weil pairing) that no
elliptic curve defined over $\Q$ satisfies the assumption that all of its 3 or
5 torsion is rational.  Therefore we must be taking elliptic curves over larger
number fields.   However, given an elliptic curve $E$ over $\Q,$ it is still
may be worthwhile to `base-change' to a number field $K$ where the assumption
does hold; performing descent over $K$ may determine, for example,  the
Mordell-Weil rank of $E$ over $K,$ thus bounding the Mordell-Weil rank of $E$
over $\Q.$

Work is currently in progress, involving the author and others, to explicitly
perform 3-descent on $E$ (in Weierstrass form) without any assumptions of
rationality of the 3-torsion points of $E.$

The paper is organised as follows:  first we will give the general set-up,
reviewing the mechanics of $p$-descent on elliptic curves and our approach,
which will consist of two basic parts.  Next we will work through the cases
$p=3$ and $p=5$ separately.

\section{The General Set-up}

Fix a number field $K$ and an elliptic curve $(E, \oh_E \in E(K))$ defined over
$K,$ where $\oh_E \in E(K)$ is the origin of the group law for $E(K).$  Then
the Mordell-Weil group $E(K)$ is known to be a finitely generated abelian
group, so we can write it as a product of a free part and a torsion part:
$E(K)= \Z^r \times E(K)_{tors}.$   For a prime $p\in \N,$ we then know that 
$E(K)/ p E(K)= (\Z/p \Z)^r \times E[p](K)$ is finite, and its dimension as an
$\F_p$-vector space is $r+\epsilon,$ where  $\epsilon= 0, 1$ or 2.  A major and
historical goal of number theorists is to compute $r.$  The theory of
$p$-descent is one approach.

Indeed, `(The first) $p$-descent on $E$' can be described as an attempt to
compute the finite group $Sel(E)[p];$  this group sits in a commutative diagram
whose rows are exact sequences and where each group on the top row is a
subgroup of the corresponding group on the bottom row:   
$$\begin{array}{ccccccccc} 1 & \rightarrow & E(K)/ p E(K) & \rightarrow &
Sel(E)[p] & \rightarrow & \Sha(E)[p]  & \rightarrow & 1\\ & & = \downarrow & &
\downarrow & & \downarrow & &\\ 1 & \rightarrow & E(K)/ p E(K) & \rightarrow &
H^1(K, E[p]) & \rightarrow & H^1(K, E) & \rightarrow & 1\\ \end{array},$$

The elements of the torsion cohomological group $H^1(K, E)$ have long been
(abstractly) identified with (isomorphism classes of) principal homogeneous
spaces for $E.$  A nontrivial class of $H^1(K, E)$ is represented by a genus
one curve with no points defined over $K$ and whose Jacobian elliptic curve is
isomorphic to $E.$  The right top group $\Sha(E)$ is a subgroup of $H^1(K, E)$
and as such, an element of $\Sha(E)$ is characterized as a class
represented by a curve $C$ for which, for every place $v$ of $K,$ there {\it
is} a $K_v$-rational point of $C.$  A famous conjecture states that $\Sha(E)$
is finite.

If in fact $\Sha(E)[p]$ is trivial, then 'knowing $Sel(E)[p]$' is equivalent to
'knowing $E(K)/ p E(K),$' which by the above discussion is very close to
computing $r;$ indeed if one has a handle on the $p$-torsion of $E,$ then it 
is the same.  If on the other hand $\Sha(E)[p]$ is non-trivial, then one
proceeds to the second level of descent involving the following commutative
diagram: $$\begin{array}{ccccccccc} 1 & \rightarrow & E(K)/ p^2 E(K) &
\rightarrow & Sel(E)[p^2] & \rightarrow & \Sha(E)[p^2]  & \rightarrow & 1\\ & &
 \downarrow & & \downarrow & & \downarrow & &\\ 1 & \rightarrow & E(K)/ p E(K)
& \rightarrow & Sel(E)[p] & \rightarrow & \Sha(E)[p] & \rightarrow & 1\\
\end{array}.$$ The key point is that the right-most vertical arrow is
`multiplication by $p$,' so if $\Sha[p^2]=\Sha[p],$ then the image of the
middle vertical map all comes from the left.  In other words, we can think of
the initial contribution from $\Sha[p]$ as noise which we've filtered out by going
to the second level of descent. Notice if $\Sha(E)$ is finite, descent will
eventually be effective.

Since we've already understood 2-descent, it might seem natural to next try to
understand 4-descent.  However, another approach to `filtering out' the noise
of $\Sha(E)[2]$ is to switch {\em primes}; that is, we hope that if
$\Sha(E)[2]$ is non-trivial, than $\Sha(E)[3]$ will be trivial.  Moreover, the
actual computations grow very quickly;  we therefore save ourselves work by
restricting to $p=3$.

A more refined description of the goal of $p$-descent is then to  not only
compute the middle group $Sel(E)[p]$ but to piece together what part `comes
from the left' and what part `contributes to the right;'  thus descent breaks
naturally into two basic parts. The first part is simply an explicit
description of the first map.  The second part is an explicit identification of
elements of $Sel(E)[p]$ with smooth models of the appropriate genus one curves
inside projective space.  The reason we want explicit models is that we can
then actually look for `local points,' that is, points over the field $K_v,$
and thus determine whether that curve contributes to the right.

Now assume $K$ is large enough so that $E[p](K) = E[p](\Kbar),$ and choose an $\F_p$
basis $S, T$ of $E[p](K)$ with fixed Weil pairing. Then the Galois cohomology group
$H^1(K, E[p])$ can be identified with the group $K^*/K^{*p} \times K^*/K^{*p}$ (see
Lemma 3.1 on page 4 of \cite{me2}), and its subgroup $Sel(E)[p]$ can be identified
with a finite subgroup of $K^*/K^{*p} \times K^*/K^{*p}.$

For the first part of descent (see Chapter X of (\cite{Sil})), we find a pair of
rational functions  $(f_S, f_T)$ on $E$  which, when evaluated at a point of $E(K),$
give its image in $H^1(G, E[p]) \cong K^*/K^{*p} \times K^*/K^{*p}.$  Note this will
depend on the chosen basis above. By Theorem 1.1d on page 278 of \cite{Sil}, it is
sufficient for  the functions $f_S$ and $f_T$ to satisfy $div(f_S) = p \cdot (S) - p
\cdot (O_E)$ and $div(f_T) = p \cdot (T) - p \cdot (O_E)$ respectively;  moreover,
they can be chosen to satisfy $f_S \circ [p] = g_S^p$ and $f_T \circ [p] = g_T^p$
for some rational functions $g_S$ and $g_T.$

\begin{lemma}\label{expansionato} The expansions of $f_S$ and $f_T$  with respect to
a local parameter at $O_E$ can be chosen to  have leading coefficients which are
perfect $p$th powers, i.e. so $f_S = \frac{a}{t^{p}} + \dots$ and  $f_T =
\frac{b}{t^p} + \dots$ where both $a$ and $b$ are perfect $p$th powers and where $t
\in \oh_{E, O_E}$ is  a parameter in the local ring at $O_E.$  Here the `$\dots$'
refer to `higher order terms.' \end{lemma}

\bpf First we prove that locally the expansions look like: $$f_S= \frac{a}{t^{p}} +
\dots, \hspace{.2in} g_S = \frac{c}{t} + \dots, \hspace{.2in} \mbox{and} 
\hspace{.2in} t \circ [p] = d \cdot t +\dots.$$  The first is because we know $f_S$
has a pole of order $p$ at $O_E,$ the second because $g_S$ has a simple pole at
$O_E;$ in fact  $div(g_S) = \sum_{p P_i = S} (P_i) -\sum_{p Q_i = O} (Q_i),$ and we
get only one copy of $O_E$ on the right.  Finally, $t \circ [p]$ has a simple zero
at $O_E$ since $[p] O_E = O_E$ and $[p]$ is \'etale when char$(K) \not = p.$ We know
that  $f_T \circ [p] = g_T^p,$ and a quick calculation shows this is equivalent to
$a$ being a perfect $p$th power. \epf

The background for the second part is to be found in \cite{me2},
and is summarized as follows.

We want to find an efficient model of $C.$  However, there isn't always a degree
three line bundle on $C,$ so we can't always model $C$ as a smooth cubic in
$\PP^2.$  The order of $C$ as an element of $H^1(K, E)$ (called the `period' of $C$)
is closely related to the smallest degree of a line bundle $\LL$ on $C$ (called the
`index' of $C$).  To explain this better, we will view the Selmer group $Sel(E)[p]$
as a subgroup of the cohomology group $H^1(K, E[n]).$ By Proposition 2.2 on page 3
of \cite{me2},  elements of  $H^1(K, E[n])$ correspond to diagrams $C \rightarrow
S$  where $C$ is a period $n$ principal homogeneous spaces of the elliptic curve $E$
as above, and where $S$ is a Brauer-Severi variety of dimension $n-1.$  These
diagrams are twists of a fixed diagram $E \rightarrow \PP^{n-1},$ given by the
divisor $n \cdot O_E.$  Under the natural map from $H^1(K, E[p])$ to $H^1(K, E),$
the diagram $C \rightarrow S$ goes to $C.$  The other forgetful map, sending $C
\rightarrow S$ to $S,$ is a quadratic map from $H^1(K, E[p])$ to 
$H^1(K,\pgl_p)\subset H^2(K,  \GG_m)[p]=Br(K)[p]$ (for this inclusion see p. 158
of \cite{Serre}) and is called the \emph{period-index obstruction map}. The
obstruction is trivial exactly when $S$ is isomorphic over $K$ to $\PP^{p-1},$ and
in this case we say that the diagram $C \rightarrow S$ has its period equal to its
index.  In our situation we will make use of our identification $H^1(K, E[p]) \cong
K^*/K^{*p} \times K^*/K^{*p};$ then the period-index obstruction for an element
$(a,b) \in K^*/K^{*p} \times K^*/K^{*p}$ is (Proposition 3.4 on page 6 of
\cite{me2}) its associated `Hilbert symbol' $(a, b)_{Hilb, p}.$ In particular, $(a,
b)_{Hilb, p}$ is trivial exactly when  $b$ is in the image of the norm map from the
field $K(\alpha)$ to $K$, where $\alpha$ is chosen such that $\alpha^p =a.$

A crucial fact that we will take advantage of is that
the elements of $Sel(E)[p],$ which a priori correspond to diagrams
$C \rightarrow S,$ always have $S \cong_K \PP^{p-1};$  in other
words, diagrams $C\rightarrow S$ where $C \in \Sha(E)$
{\it always} have trivial period-index obstruction
(for a proof see the Remark on page 3 of \cite{me2}).

In this article we will perform the second part of descent not only for elements of
$Sel(E)[p]$ but for any element of $H^1(K, E[p])$ which has trivial period-index
obstruction.  Indeed we will start with a  ``trivialisation'' of a Hilbert symbol
(the data of an element $\beta \in K[x]/(x^p-a)$ whose norm is $a$), and produce a
model for the corresponding genus one curve under the isomorphism $H^1(K, E[p])
\cong K^*/K^{*p} \times K^*/K^{*p}.$

\section{Explicit Descent over $X(3)$}
\subsection{The image of $E(K)/3 E(K)$}

Let $K$ be a field whose characteristic is different
from 3 and which contains a primitive 3rd root of unity $\zeta.$  
Let $E$ be an elliptic curve over $K$ with
full 3-torsion given by the following equation:
$$E: X^3+Y^3+Z^3+ \lambda X Y Z = 0.$$
Define the origin $O_E$ to be the point $(1; -1 ; 0)$
and a basis of 3-torsion points $<S, T>$  to be
$S=(1; \zeta; 0), T=(1; 0 ; -1),$ where
$\zeta \in K$ is a fixed cube root of unity.

By results in \cite{me}, there exist matrices in $\pgl_3$
which act as ``translation by $S$ and $T$" on $E$ and which 
we will denote by $M_S$ and $M_T.$  We leave it to the reader 
to prove these matrices are $D_3$ and $M_{1,3}$ respectively, where
$D_n$ is the diagonal matrix $diag(1, \zeta_n, \zeta_n^2, \dots, \zeta_n^{n-1})$
for some primitive $n$ root of unity $\zeta_n,$ (here take $\zeta_3= \zeta$)
and where 
$$M_{a, n}= \left( \begin{array}{ccccc}
                 0 & 1 & 0 & \dots & 0\\
                 0 & 0 & 1 & \dots & 0\\
		 \vdots & \vdots & \vdots &  & \vdots \\  
		 0 & 0 & 0 & \dots & 1\\      
                 a & 0 & 0 & \dots & 0
                  \end{array} \right).$$\label{matrices}

We will find a pair of rational functions $(f_S, f_T)$ on $E$ 
which when evaluated at a point of $E(K)$ gives its image in
$H^1(G, E[3]) \cong K^*/K^{*3} \times K^*/K^{*3}.$
By Lemma \ref{expansionato}, the function $f_S$ needs to satisfy
the condition that $div(f_S) = 3 \cdot (S) - 3 \cdot (O_E);$
this is satisfied
trivially by the quotient of the hyperplane at $S$ by the hyperplane at 
$O_E,$
since they are both flex points.  The hyperplane at $O_E$ is easily seen
to be given by $3X+3Y-\lambda Z,$ and the other hyperplanes are in fact
translates of this one by the above matrices; this observation doesn't
save much time in $\PP^2,$ since hyperplanes are easy to compute, but will in
the $n=5$ case.  Similarly $f_T$ is, up to a constant,
a ratio of hyperplanes, and we get
$f_S= \kappa_S \frac{3 \zeta^2 X+3 \zeta \, Y- \lambda \, Z}
{3X+3Y-\lambda Z}$ and $f_T=\kappa_T \frac{3 X- \lambda Y
+ 3 Z}{3X+3Y- \lambda Z}$ for some constants $\kappa_S$ and $\kappa_T.$

\begin{proposition}
$f_S= (\lambda^3 + 27) \frac{3 \zeta^2 X+3 \zeta \, Y- \lambda \, Z}
{3X+3Y-\lambda Z}$ and $f_T=(\lambda^2 - 3 \lambda +9) \frac{3 X- \lambda Y
+ 3 Z}{3X+3Y- \lambda Z}$. 
\end{proposition}

\bpf
By Lemma \ref{expansionato},
we want then to compute the expansion of $\frac{3 \zeta^2 X+3 \zeta \, Y- 
\lambda \, Z}{3X+3Y-\lambda Z}$ at the origin $O_E=(1; -1 ; 0).$
Here we can work in affine coordinates by setting $X=1,$ since this
is true locally.  Then $F: 1+ Y^3+Z^3+ \lambda YZ =0.$
Moreover, since $O_E$ takes on a non-zero value
at the hyperplane at $S,$ we can just evaluate there to get
$(3 \zeta^2 X+3 \zeta \, Y- \lambda \, Z)|_{O_E}= 3 \zeta^2-3\zeta.$
Note that $\zeta^2 - \zeta = \sqrt{-3}$, so 
$3 \zeta^2-3\zeta = (-\sqrt{-3})^3$
is a cube in $K.$  So in fact we can completely ignore this term. 
We're left with $\frac{1}{3X+3Y-\lambda Z}$ which has a triple zero at $O_E.$
In other words
$$3+3Y- \lambda Z \cong a Z^3 + \dots,$$
where $\cong$ signifies that we are working in  $\oh_{E, O_E}.$
Multiply the above by the function $Y,$ a nonzero funtion at $O_E,$ to get
$3 Y+3Y^2- \lambda YZ \cong a Y Z^3 + \dots;$
since $-\lambda Y Z \cong 1 + Y^3 + Z^3$ (we're working modulo $F$)
we substitute to get
$$1+3Y + 3Y^2 + Y^3 + Z^3 = (1+Y)^3 + Z^3 \cong a Y Z^3 + \dots.$$
The function $1+Y$ has a zero at $O_E,$ so it's an alternative parameter for the
local ring $\oh_{E, O_E}:$ $1+Y \cong b Z +\dots.$
but we already have the tangent line equation which tells us that
$3(1+Y) \cong \lambda Z + \dots,$ i.e. $b=\lambda/3.$
Replace $(1+Y)^3$ above now by $(\lambda/3) Z^3:$
$$Z^3 ( 1+ (\lambda/3)^3) \cong a Y Z^3.$$ Evaluate at $Y=-1$ to get
$a = -(1+(\lambda/3)^3).$  Our original constant then is $1/a,$ which
modulo cubes is seen to be $\frac{1}{27+\lambda^3}.$  Finally, to normalize
we want the function $f_S$ to have a leading coefficient which is 1.
In other words, we need to multiply the ratio of the two hyperplanes by
the constant $27 + \lambda^3.$
For $f_T,$ we are already almost done- the only difference is the value
of the numerator at the origin:
$(3 X- \lambda Y + 3 Z)|_{O_E}= 3 + \lambda.$  The leading coefficient then
is $\frac{3+\lambda}{27+\lambda^3} = \frac{1}{\lambda^2 - 3 \lambda + 9}.$
\epf

\subsection{Models of Genus One curves in $\PP^2$}

Let $E$ be an elliptic curve over the field $K.$  Assume $char(K) \not = 3$ 
and that $E[3](K) = E[3](\Kbar).$  Choose a basis $\langle S, T \rangle$ 
of $E[3]$ and identify $H^1(K, E[3])$ with $K^*/K^{*3} \times K^*/K^{*3}.$
Take $(a, b) \in H^1(K, E[3])$ whose corresponding Hilbert symbol is trivial.
If both $a$ and $b$ are cubes, the element $(a,b)$ represents $E.$
We will first assume that $a$ is not a cube in $K;$
the other case will be dealt with at the end of this section.
Define $\alpha, \beta \in \Kbar$ such that
$\alpha^3 = a$ and $\beta = \sum_{i=0}^2 \beta_i \alpha^i$ is
such that $\N_{K(\alpha)/K} (\beta) = b.$  Moreover define
$\sigma \in G(K(\alpha)/K)$ to be the generator of that cyclic Galois group
such that $\sigma(\alpha) = \alpha \zeta.$ 
Let $Tr$ be the trace function of $K(\alpha)$ to $K,$
and let $u= \beta/{\sigma^2 ( \beta)}.$

\begin{theorem}\label{p2formula}
The curve $(a,b)$ is given by
$$(Tr (u) + \lambda) \, ( a^2 X^3 + a Y^3 + Z^3) + 
3 \, Tr(\alpha u ) \, ( a X^2 Z + a Y^2 X + Z^2 Y) \, +$$  
$$3 \, Tr(\alpha^2 u) \, (a X^2 Y + Y^2 Z + Z^2 X) \, + 
3 \, ( 2a \, Tr(u) - \lambda a) \, X Y Z =0 .$$
\end{theorem}

\bpf
By results in \cite{me}, the action of the 3-torsion points
of $E$ on a a genus one curve $C$ embedded in $\PP^2_K$ can be
represented as automorphisms of $\PP^2_K,$ in other words
as elements of $\pgl_3(K).$  Moreover, if $C$ is represented in
$H^1(G, E[3])$ by the pair $(a, b)$ which depends on a chosen
basis $<S,T>$ then the determinants of the matrices representing
``translation by $S$'' and ``translation by $T$'' are $a$ and
$b$ respectively.  Moreover the Weil pairing is given by
the commutator of lifts to $\gl_3(K).$  With that in mind,
we will search for $C \leftrightarrow (a,b)$ by finding
cubics which are invariant under the action of standard matrices 
representing this 3-torsion action.
Define $M_S$ to be the matrix $M_{a, 3}$
and $M_T$ to be the matrix $D_3 \, [\beta_0 I+
\beta_1 M_S + \beta_2 M_S^2],$  if $\beta= \beta_0 + \beta_1 \alpha
+ \beta_2 \alpha^2.$
The determinant of $M_S$ is $a,$ the determinant of $M_T$
is $b,$ and the commutator $[M_S, M_T]$ is $\zeta I.$
A cubic which is invariant under the action of $M_S$ 
but with no fixed points must be of the form
$$F= A \, ( a^2 X^3 + a Y^3 + Z^3) + 
B \, ( a X^2 Z + a Y^2 X + Z^2 Y) \, +$$  
$$C \, (a X^2 Y + Y^2 Z + Z^2 X) \, + 
3 \, D \, X Y Z =0 .$$  This is because $M_S$ acts linearly on the 
10-dimensional space of cubics.  There are 3 eigenspaces of dimensions
3, 3, and 4, and the first two have zeroes at the fixed points of $M_S$,
whereas the last eigenspace does not and is generated by the above four cubics.

On the other hand we also insist that $F$ be invariant under the action of 
$M_T.$  To ease computations we introduce the following notation: 
fix eigenvectors 
$v_i = (1, \alpha \zeta^i, \alpha^2 \zeta^{2\, i})$ of $M_S.$
Then $M_S v_i= \alpha \zeta^i v_i.$ 
The four coefficients $A, B, C,$ and $D$ are linear combinations
of $F(v_0), F(v_1), F(v_2),$ and $\TT(v_0, v_1, v_2),$
where $\TT$ is the trilinear form associated to $F,$ as follows:
$$\begin{array}{cccc}
\left( \begin{array}{c} F(v_0) \\ F(v_1)\\ F(v_2) \\ \TT(v_0, v_1, v_2)
\end{array} \right) & = &  \left(
\begin{array}{cccc} 3a^2 & 3 a^{5/3} & 3 a^{4/3} & 3a \\
3a^2 & 3 a^{5/3} \zeta^2  & 3 a^{4/3} \zeta & 3a \\ 
3a^2 & 3 a^{5/3} \zeta & 3 a^{4/3} \zeta^2 & 3a \\
18 a^2 & 0 & 0 & -9a \end{array} \right) &
\left( \begin{array}{c} A \\ B \\ C \\ D
\end{array} \right) \end{array}.$$
Now it is easy to see how $M_T$ acts on the $F(v_i):$
$$F^{M_T}(v_i)=F(M_T v_i) = F(D \, [\beta_0 I+
\beta_1 M_S + \beta_2 M_S^2] v_i)=$$  
$$F(D \, (\beta_0 + \beta_1 \alpha \zeta^i + \beta_2 \alpha^2 \zeta^{2 \, i}) v_i)=
F(\sigma^i (\beta) D v_i) = \sigma^i (\beta)^3 F(v_{i+1}).$$
Similarly, $$\TT^{M_T}(v_0, v_1, v_2) = \TT(\beta v_1, \sigma(\beta) v_2,
\sigma^2(\beta) v_0) = b \, \TT(v_0, v_1, v_2).$$
The fact that $F$ is invariant under the action of $M_T$
is equivalent to the projective point
$P=(F(v_0); F(v_1); F(v_2); \TT(v_0, v_1, v_2))$ being fixed by
$M_T.$  We have seen that $M_T(P) =
(\beta^3 F(v_1); \sigma(\beta)^3 F(v_2); \sigma^2(\beta)^3 F(v_0);
b \, \TT(v_0, v_1, v_2)).$
For some $\mu \neq 0$ we have 
$F(v_1) = \frac{\mu}{\beta^3} F(v_0)$ and $F(v_2)= \frac{\sigma^2(\beta)^3}{\mu} F(v_0).$
Moreover, the Jacobian of the above curve is (see \cite{me}):
$$X^3 + Y^3 + \prod_{i=0}^2 F(v_i) Z^3  +
\TT(v_0, v_1, v_2) XYZ=0,$$ in other words
$$X^3 + Y^3 + \frac{F(v_0)^3 \sigma^2(\beta)^3}{\beta^3} Z^3  +
\TT(v_0, v_1, v_2) XYZ=0;$$ setting $F(v_0) = \frac{\beta}{\sigma^2(\beta)}$ 
(note that $F(v_0) \neq 0$ because $F$ has no fixed points under the action of $M_S$)
and renaming $\TT(v_0, v_1, v_2) = \lambda,$ the Jacobian is exactly $E.$
Note that $\sigma(F(v_i)) = F(\sigma(v_i)) = F(v_{i+1}),$ so
$F(v_1) =  \frac{\sigma(\beta)}{\beta}$ and  $F(v_2) =  
\frac{\sigma^2(\beta)}{\sigma(\beta)}.$
To finish the proof, we need to invert the above matrix to find
the coefficients $A, B, C,$ and $D$ in terms of the $F(v_i)'s$ and
$\TT:$
$$\begin{array}{cccc}
\left( \begin{array}{c} A \\ B \\ C \\  D 
\end{array} \right) & = &
\frac{1}{27 a^2}
\left(\begin{array}{cccc}
1 & 1 & 1 & 1 \\
3 \alpha & 3 \alpha \zeta & 3 \alpha \zeta^2 & 0 \\
3 \alpha^2 & 3 \alpha^2 \zeta^2 & 3 \alpha^2 \zeta & 0 \\
2 a & 2 a & 2 a & -a 
\end{array}\right) & 
\left( \begin{array}{c} F(v_0) \\ F(v_1) \\ F(v_2) \\ \TT(v_0, v_1, v_2) 
\end{array} \right)
\end{array}.$$

We can ignore the factor of $\frac{1}{27 a^2},$ since we're working 
projectively.  Then $A= F(v_0) + F(v_1) + F(v_2) + \TT(v_0, v_1, v_2) = 
\frac{\beta}{\sigma^2(\beta)} + \frac{\sigma(\beta)}{\beta}
+ \frac{\sigma^2(\beta)}{\beta} + \lambda = Tr(u)+ \lambda,$ 
and similarly for the other coefficients. 
\epf

\noindent{\bf When $a$ is a cube}
The formula in Theorem \ref{p2formula} was computed assuming $a$ 
is not a cube in the base field $K,$  but it holds in
the case that $a$ is a cube as well.  We can view the formula as
holding over the algebra $K[\alpha]/(\alpha^3 -a),$ which is a field if
$a$ is not a cube; in other words we view $\alpha$ as strictly a symbol.
The norm function on this algebra
takes an element $s =s_0  +s_1 \alpha  +s_2 \alpha^2 $ to the expression
$\N(s)= \prod_{i=0}^2 (s_0 +s_1 \alpha \zeta^{i} +s_2 \alpha^2 \zeta^{2i}) =
s_0^3 + s_1^3 a + s_2^3 a^2 - 3 s_0 s_1 s_2 a.$
Similarly, the trace function has
$Tr(s)= \sum_{i=0}^2 (s_0 +s_1 \alpha \zeta^{i} +s_2 \alpha^2 \zeta^{2i}) =
3 \cdot s_0,$ since $1+\zeta +\zeta^2=0.$
In other words, we can define the linear automorphism `$\sigma$' to take
$\alpha$ to $\zeta \alpha,$ and $\N(s) = s \sigma(s) \sigma^2 (s).$
In the case of representing a curve corresponding to the pair
$(a,b)$ for $a$ a cube, we may take $a= 1$ and we want to find
$\beta = \beta_0   +\beta_1 \alpha  +\beta_2 \alpha^2$ such that
$\N(\beta) = \beta \sigma(\beta) \sigma^2(\beta) =b.$  There are many
solutions:  for example, we can
take $\beta_0 = \frac{b+2}{3}$ and $\beta_1 =\beta_2 = 
\frac{b-1}{3}$ if those are nonzero quantities.  
Then as above we can formally compute $Tr(u) = (b-1)^2/ b +3,$
$Tr(\alpha u) = (b-1) (b-\zeta^2)/ b,$ and 
$Tr(\alpha^2 u) = (b-1)(b-\zeta)/ b.$

\section{Explicit Descent over $X(5)$}
\subsection{The image of $E(K)/5 E(K)$}

Let $K$ be a field whose characteristic is different
from 5 and which contains a primitive 5th root of unity $\zeta.$
Let $E= E_{\lambda}$ be the elliptic curve over $K$ 
given by the equations
$\lambda x_i^2 +\lambda^2 x_{i-2} x_{i+2} - x_{i-1} x_{i+1}$ 
for $i$ between 0 and 4 and whose origin is given by 
$O_E= ( 0; \lambda; 1; -1; -\lambda).$

\rk
The curves $E_\lambda$ for $n=3,5$ are also
considered in \cite{fisher}, where the Cassels-Tate 
pairing is used to compute $S^{(n)}(E_\lambda/\Q)$
for all non-cuspidal $\lambda \in \Q$.  
Equations for $E_\lambda$ for $n=5$ in 
Weierstrass form may be found in \cite{fisher} or \cite{Beaver}.

The curve $E$ as above has full 5-torsion over $K$ and the matrices
which act as translation by a $S$ and $T,$ for a basis $\langle S, T
\rangle,$ are given by $D_5$ and $M_{1,5}$ respectively (see notation
on page \pageref{matrices}).

We will find a pair of rational functions $(f_S, f_T)$ on $E$ 
which when evaluated at a point of $E(K)$ gives its image in
$H^1(G, E[5]) \cong K^*/K^{*5} \times K^*/K^{*5}.$
By Lemma \ref{expansionato}, the function $f_S$ can be chosen to satisfy
$div(f_S) = 5 \cdot (S) - 5 \cdot (O_E);$
then we can take $f_S$ to be a scalar multiple of the quotient of the 
hypertangent plane at $S$ by the hypertangent plane at $O_E,$
since both $S$ and $O_E$ are hyperflex points (to see this, note that
$E$ is a degree 5 curve in $\PP_K^4$ and that
the hyperplane $x_0=0$ goes through the points $i \cdot S$ for
$i$ between 0 and 4, which means that the divisor giving the
embedding $E \rightarrow \PP^4$ is linearly equivalent to
$(O_E) + (S) + (2S) + (3S) + (4S) \equiv 5 \cdot (O_E)\equiv 5 \cdot (S)).$

\begin{proposition}
The hypertangent plane at the origin of $E$ is given by
$$H_{O_E}:\alpha x_0 + \beta(x_1 + x_4) + \gamma(x_2 + x_3),$$
where  $\alpha = \lambda^{10} - 14 \lambda^5 -1, 
\beta = -5 \lambda^2 (1+2 \lambda^5),$ and $\gamma = 5\lambda^3
(\lambda^5 -2).$ 
\end{proposition}

\bpf
An easy calculation (in Maple for example) verifies that the above hyperplane 
intersects $E$ only at $O_E.$  In order to find the above,
we computed a local parameterization of the curve $E$ in the
local ring $\oh_{E, O_E}$ using the equations for $E$ and Maple.
\epf

To finish finding the equations $f_S$ and $f_T,$ we simply
find the hypertangent planes at $S$ and $T$ (these are translates 
of $H_{O_E}$ by the 5-torsion matrices $D_5$ and $M_{1,5}$) 
and evaluate them at $O_E.$   We then scale
$H_S/H_{O_E}$ by the appropriate function of $\lambda$ so that its
leading coefficient in the expansion at $O_E$ is a perfect fifth power:

\begin{proposition}
The hypertangent planes at $S$ and at $T$ are given by
$$H_S:\alpha x_0 + \beta(x_1 \zeta^4 + x_4 \zeta) + 
\gamma(x_2 \zeta^3+ x_3 \zeta^2)$$ and
$$H_T:\alpha x_4 + \beta(x_0 + x_3) + \gamma(x_1 + x_2).$$
The rational functions $f_S$ and $f_T$ are given by
$$f_S= \frac{[(\lambda^2+\lambda-1)(\lambda^4 - 3 \lambda^3 + 4 \lambda^2 - 2 
\lambda +1)(\lambda^4 + 2 \lambda^3 + 4 \lambda^2 +3 \lambda+1)]^2}
{5 \; \lambda \; (\zeta-\zeta^4) \; (\lambda^5-2)}  \frac{H_S}{H_{O_E}};$$ 
$$f_T=\frac{\lambda^2(\lambda^4 - 3 \lambda^3 + 4 \lambda^2 - 2 
\lambda +1)^2 (\lambda^4 + 2 \lambda^3 + 4 \lambda^2 +3 \lambda +1)}{
(\lambda^2+\lambda-1)}
\frac{H_T}{H_{O_E}}.$$ 
\end{proposition}

\subsection{Models of Genus One curves in $\PP^4$}

Let $E$ be an elliptic curve over the field $K.$  Assume $char(K) \not = 5$ 
and that $E[5](K) = E[5](\Kbar).$  Choose a basis $\langle S, T \rangle$ 
of $E[5]$ and identify $H^1(K, E[5])$ with $K^*/K^{*5} \times K^*/K^{*5}.$
For $(a, b) \in H^1(K, E[5])$ whose corresponding Hilbert symbol is trivial,
assume as in the previous section that $a$ is not a perfect fifth power
in $K$, and define $\alpha, \beta \in \Kbar$ such that
$\alpha^5 = a$ and $\beta = \sum_{i=0}^4 \beta_i \alpha^i$
such that $\N_{K(\alpha)/K} (\beta) = b.$  Moreover define
$\sigma \in G(K(\alpha)/K)$ to be the generator of that cyclic Galois group
such that $\sigma(\alpha) = \alpha \zeta.$ 

\begin{theorem}
There exists a unique 5-dimensional $K$-vector space of quadrics
in the variables $x_i, i = 0, \dots 4$, denoted $\V,$
satisfying the following:
\begin{itemize}
\item $\V$ defines a smooth degree 5 genus one curve $\C$ in $\PP^4.$
\item The Jacobian of $\C$ is isomorphic over $K$ to $E,$ and in fact
the diagram $\C \rightarrow \PP^4$ corresponds to the pair
$(a, b) \in H^1(G, E[n](\Kbar)).$
\item The action of $E[n]$ on $\C,$ a principal homogeneous space of
$E,$ is determined by the action of $S$ and $T$, and these are given as
restrictions of the following automorphisms of $\PP^4:$
`action by $S$'= $M_S$ = $M_{a,5}$, and
`action by $T$'= $M_T$ = $D_5 \cdot \sum_{i=0}^4 \beta_i M_S^i.$ 
\end{itemize}
\end{theorem}

\bpf
By Proposition 2.2 on page 3 of \cite{me2}, the cohomology group 
$H^1(G, E[n](\Kbar))$ parameterizes
diagrams $C \rightarrow \PP^{n-1}$ up to isomorphism, where an isomorphism
is a map defined over $K$ from $C$ to $C$ which extends to a $K$-map on
$\PP^{n-1},$ i.e. an element of $\pgl_n(K).$  Then $\V$ determines such
a diagram, and the prescribed action of $E[n]$ fixes the model:
modifying a diagram $C \rightarrow \PP^{n-1}$ 
by an element of $\pgl_n(K)$ has the effect of 
conjugating $M_S$ and $M_T.$  Since the group generated by 
$M_S$ and $M_T$ (a twist of the Heisenberg group) is its own centraliser,
there is no non-trivial modification.
\epf

Our goal of this section is to determine $\V.$

\begin{definition} 
Define, for $i= 0, \dots, 4$, $v_i$ to be points in $\A^5_{\Kbar}$ such that
$M_S v_i = \alpha \zeta^i v_i.$ In particular we see that 
$\sigma(v_i) = v_{i+1}.$
\end{definition}

\begin{definition}
Define the action of a matrix $M$ on a function $f$ to be such that
$f^M(x) = f(M x).$  Then $f^{M_1 M_2} = (f^{M_1})^{M_2}.$
\end{definition}

\begin{lemma}\label{normalizeQ}
With notation as above,
we can choose a quadric $Q \in \V$ such that $Q(v_0) = Q(v_1)$ and 
such that $\langle Q^{M_S^i} \rangle_{i = 0 \dots 4}$ forms a basis of $\V.$
Moreover, such a $Q$ is unique up to $K$-scaling.  Therefore to determine
$\V$ it suffices to determine $Q.$
\end{lemma}

\rk  We are actually taking a fixed quadric $Q$ (not defined up to
a scalar), since we need to make sense of the nonzero quantities $Q(v_0)$
and $Q(v_1).$

\bpf
We can span $\V$ by translates of any $K$-rational quadric $Q$ 
by the action of powers of $M_S$ since the eigenvalues of 
$M_S$ acting on the space of quadrics are the fifth roots of
$a^2,$ not defined over $K.$
Note that $Q(v_i) \not = 0$ for all $i$ since if so we would have 
$Q^{M_S^i}(v_i) = Q(M_S^i v_i) = 0$ as well, in other words we would have
a point on $C$ fixed by $M_S,$ namely the projectivization of $v_i.$
Next, note that $Q^{M_S^i}(v_0) = Q(M_S^i v_0) = Q( \alpha^i v_0)
=\alpha^{2i} Q(v_0).$  Therefore if we replace $Q$ by the quadric
$Q' = \sum_{i=0}^4 a_i Q^{M_S^i},$ for $a_i \in K,$ then 
$Q'(v_0) = \sum_{i=0}^4 a_i Q^{M_S^i}(v_0) = 
\left(\sum_{i=0}^4 a_i \alpha^{2i} \right) Q(v_0),$ and likewise
$Q'(v_1) = \sum_{i=0}^4 a_i Q^{M_S^i}(v_1) = 
\left(\sum_{i=0}^4 a_i \alpha^{2i} \zeta^{2 i} \right) Q(v_1).$
So in order to have $Q'(v_0) = Q'(v_1),$ we need to find 
$a = \sum_{i=0}^4 a_i \alpha^{2i}$
so that $Q(v_0)/Q(v_1) = \sigma(a)/a.$  This is possible by Hilbert's
Theorem 90, since $Q(v_0)/Q(v_1) = Q(v_0)/\sigma(Q(v_0))$ is in the 
kernel of the norm map.  Next, $Q'$ and its
translates under $M_S$ also generate $\V.$  The $K$-rational
quadric $Q'$ is nonzero
since its value at $v_0$ is nonzero;  by the above comment its 
translates by powers of $M_S$ span $\V.$  Finally, such a $Q'$ is
unique up to an element of $K,$ since if we had both $Q$ and $Q'$
such that $Q(v_0)/Q(v_1) = Q'(v_0)/Q'(v_1) = 1,$ we could write
$Q' = \sum_{i=0}^4 a_i Q^{M_S^i}$ for $a_i \in K$ to get
$Q'(v_0)/Q'(v_1) = a/\sigma(a) \cdot Q(v_0)/Q(v_1),$  i.e.
we would have $a = \sigma(a),$ or in other words $a \in K.$
\epf

Our newly defined goal is to find $Q$ as in Lemma \ref{normalizeQ}.
We will choose a matrix $M$ by which to modify $Q$ as in
Lemma \ref{normalizeQ} so that the coefficients of $Q^M$ are easy
to manipulate.
Define $M=(v_0 v_1 v_2 v_3 v_4),$ a $5 \times 5$ matrix defined
over $K(\alpha).$  Then
$$Q^M(x) = Q(M x) = Q \left( \sum_{i=0}^4 v_i x_i \right) =
\sum_{i=0}^4 Q(v_i) \, x_i^2 + \sum_{0 \leq i < j \leq 4} B(v_i, v_j) \, 
x_i x_j,$$ where $B(w, v) = Q(w+v) - Q(w) -Q(v).$

Define $M_S' = M^{-1} M_S M$ and $M_T' = M^{-1} M_T M.$  Then we have
$(Q^M)^{M_S'} = (Q^{M_S})^M$ and $(Q^M)^{M_T'} = (Q^{M_T})^M.$
A calculation shows us 
$M_S'= \alpha D_5$ and $M_T' = M_{1,5}^{-1} \cdot 
diag(\beta, \sigma(\beta),\sigma^2(\beta),
\sigma^3(\beta), \sigma^4(\beta)).$  

Next, note that by assumption $M_T$ fixes
$\V.$  Moreover, since $M_T$ is defined over $K$, the image of $Q$ under
the action by $M_T$ is again defined over $K.$  Thus there exist 
numbers $\gamma_i \in K$ such that $Q^{M_T} = \sum_{i=0}^4 \gamma_i Q^{M_S^i}.$
Acting on that equation by $M$ we get (here let $Q' = Q^M):$
$${Q'}^{M_T'} =\sum_{i=0}^4 \gamma_i {Q'}^{(M_S') ^i}.$$  Since we know
the $M_T'$ and $M_S'$ explicitly, we can compute the left and right sides of this
equation and compare them.
$${Q'}^{M_T'}(x) =
\sum_{i=0}^4 Q(v_i) \, \sigma^{i-1}(\beta)^2 \, x_{i-1}^2 + 
\sum_{0 \leq i < j \leq 4} B(v_i, v_j) \, \sigma^{i-1}(\beta)\sigma^{j-1}(\beta)
\, x_{i-1} x_{j-1}$$ and
$${Q'}^{(M_S')^i}(x) =\alpha^{2 i} \left(
\sum_{j=0}^4 Q(v_j) \, \zeta^{2 j i} \, x_j^2 + 
\sum_{0 \leq j < k \leq 4} B(v_j, v_k) \, \zeta^{j i + k i} \, x_j x_k 
\right).$$ 

\begin{proposition}\label{formofBs}
For every $i = 0, \dots, 4$, $Q(v_i) = Q(v_0),$
$B(v_{i+1}, v_{i-1}) = 
B(v_1, v_4) \cdot \prod_{j=1}^i \sigma^{j-1} 
\left( \frac{\beta^2}{\sigma(\beta)\sigma^4(\beta)} \right),$ and
$B(v_{i+2},v_{i-2})
=B(v_2, v_3) \cdot \prod_{j=1}^i \sigma^{j-1} 
\left( \frac{\beta^2}{\sigma^2(\beta)\sigma^3(\beta)} \right).$
\end{proposition}

\bpf
We compare the coefficient of $x_0^2$ on both sides of the above equation.
On the left, we get $Q(v_1) \beta^2.$ On the right we get
$\sum_{i=0}^4 \gamma_i \alpha^{2i} Q(v_0).$  Define
$\gamma = \sum_{i=0}^4 \gamma_i \alpha^{2i},$ then we have
$Q(v_1) \beta^2 = \gamma Q(v_0).$  On the other hand we have
chosen $Q$ as in Lemma \ref{normalizeQ} so that $Q(v_0) = Q(v_1) \not = 0.$
Therefore $\gamma = \beta^2.$  We can determine the rest of the $Q(v_i)'s$ 
now since $1 = \sigma(\frac{Q(v_1)}{Q(v_0)}) = \frac{Q(v_2)}{Q(v_1)}$ etc. so
all of the $Q(v_i)$'s are equal to $Q(v_0).$

Next, compare the coefficients of $x_1 x_4:$ on the left we get
$B(v_2, v_0) \sigma(\beta) \sigma^4(\beta)$ and on the right we get
$\sum_{i=0}^4 \gamma_i \alpha^{2i} B(v_1, v_4).$  In other words we have
$\frac{B(v_2, v_0)}{B(v_1, v_4)}
= \frac{\beta^2}{\sigma(\beta) \sigma^4(\beta)}.$  Acting by $\sigma$
on both sides gives $\frac{B(v_3, v_1)}{B(v_2, v_0)}
= \sigma \left(\frac{\beta^2}{\sigma(\beta) \sigma^4(\beta)} \right),$
so $\frac{B(v_3, v_1)}{B(v_1, v_4)} = \frac{\beta^2}{\sigma(\beta) \sigma^4(\beta)}
\cdot \sigma \left(\frac{\beta^2}{\sigma(\beta) \sigma^4(\beta)} \right).$
We continue in this way.
Similarly, comparing coefficients of $x_2 x_3$ on both sides we get
$\frac{B(v_3, v_4)}{B(v_2, v_3)} = \frac{\beta^2}{\sigma^2(\beta)
\sigma^3(\beta)}$ and we finish by acting on both sides by $\sigma$
and solving for $\frac{B(v_{i+2}, v_{i-2})}{B(v_2, v_3)}.$
\epf 

Now we have only to determine the three unknowns
$Q(v_0), B(v_1, v_4), $ and $B(v_2, v_3);$ moreover, since we
are actually working projectively, we only need to know two of them,
or more precisely it is adequate to know the ratios
$\frac{Q(v_0)}{B(v_1, v_4)}$ and $\frac{Q(v_0)}{B(v_2, v_3)}.$ 
Define $E_{\lambda}$ to be the elliptic curve given by the equations
$\lambda x_i^2 +\lambda^2 x_{i-2} x_{i+2} - x_{i-1} x_{i+1}$ 
for $i$ between 0 and 4 and whose origin is given by 
$\oh_E= ( \lambda; -1; 1; -\lambda; 0).$

\begin{theorem}
The Jacobian of $\C$ is $E_{\lambda},$ where 
$$\lambda = -\frac{Q(v_0)}{B(v_2, v_3)} 
\frac{\sigma^3(\beta)\sigma^4(\beta)}{\beta \sigma(\beta)}.$$
\end{theorem}

\bpf
By \cite{me}, Theorem 4.2 on page 37, the Jacobian of $\C$ is given
as $E_A$ given by the quadrics
$S_0=x_0^2-x_2 x_3 + x_1 x_4, S_1= x_1^2-x_0 x_2 +A x_3 x_4,
S_2= x_2^2 -x_1 x_3 -A x_0 x_4, S_3= x_3^2 - x_0 x_1 - x_2 x_4,$
and $S_4 = A x_4^2 + x_1 x_2 - x_0 x_3$ \footnote{This corrects
a minus sign error in that paper, namely the coefficient of $x_4^2$ in $S_4.$}
for a parameter 
$$A = \prod_{i=0}^{4}\frac{Q(v_i)}{B(v_i, v_{i+1})}.$$
A calculation using Proposition \ref{formofBs} shows that the above
$A$ in this case is the fifth power of  $\frac{Q(v_0)}{B(v_2, v_3)} 
\frac{\sigma^3(\beta)\sigma^4(\beta)}{\beta \sigma(\beta)}.$
Finally, When $A$ is a perfect fifth power, say of $-\lambda,$ the
$K$-rational map $diag(1, -\lambda, \lambda, -1, \lambda^{-2})$
maps $E_A$ to $E_{\lambda}.$
\epf

\begin{lemma}
$$\frac{Q(v_0)^2}{B(v_1, v_4) B(v_2, v_3)} = - 
\frac{ \beta^2 \sigma (\beta)}{\sigma^3(\beta) \sigma^4 (\beta^2)}.$$
\end{lemma}

\bpf  This follows from Lemma 4.3 on page 139 of \cite{me} and
using Proposition \ref{formofBs}.  Loosely speaking, this is a 
condition on the intersection of five quadrics to form
a smooth genus one curve; note that five quadrics in general position
do not intersect.
\epf

\begin{corollary}
$$\frac{B(v_1, v_4)}{Q(v_0)} = \lambda \cdot \frac{\sigma^4(\beta)}{\beta}.$$
\end{corollary}

We now have all of the coefficients of $Q'$ in terms of $Q(v_0).$  We
can divide out by this non-zero term to get

\begin{proposition}
Let $u_1 = \frac{\sigma^4(\beta)}{\beta}$ and
$u_2 = \frac{\sigma^3(\beta)\sigma^4(\beta)}{\beta \sigma(\beta)}.$
Then $Q'$ is given by:
$$\sum_{i=0}^4 x_i^2 + \lambda \sum_{i=0}^4 x_i x_{i+2} \; 
\cdot \; \sigma^{i+1} u_1
- \lambda^{-1} \sum_{i=0}^4 x_i x_{i+1} \; \cdot \; \sigma^{i+3} u_2.$$
\end{proposition}

Finally, since $Q'=Q^{M},$ we can recover $Q$ as ${Q'}^{M^{-1}}.$
A calculation shows:

\begin{theorem}  Write
$$Q = \sum_{0 \leq i \leq j \leq 4} a_{i \; j} \; x_i x_j.$$ 
Then the $a_{i j}$ are given as follows.  For $i=j$ we have
$$a_{i \; i} = Tr \left( \frac{1}{\alpha^{2i}} \right) + 
\lambda \; Tr \left( \frac{u_1}{\alpha^{2i}} \right)  - \lambda^{-1}
\; Tr \left( \frac{u_2}{\alpha^{2i}} \right) $$ and for $i \not = j$ we have
$$a_{i \; j} = Tr \left( \frac{2}{\alpha^{i+ j}} \right)  + 
\lambda \; Tr \left( (\zeta^{i-j} + \zeta^{j-i} )\frac{u_1}{\alpha^{i+j}}
 \right) - \lambda^{-1} \;
Tr \left( (\zeta^{2i-2j} + \zeta^{2j-2i} )\frac{u_2}{\alpha^{i+j}} \right) .$$
\end{theorem}

\rk  The question of solving ``norm equations,'' that is, developing an 
algorithm to find $\beta$ as in the above theorem, has been studied extensively by Denis 
Simon \cite{simon2}, written up in the new book by Henri Cohen \cite{Cohen}, 
and implemented in Pari.

\end{document}